\documentclass[11pt, aos, authoryear]{article}
\RequirePackage[aos,amsthm,amsmath,natbib]{imsart}
\RequirePackage{hypernat}
\usepackage{amssymb, epsf, here,amsfonts,latexsym}
\usepackage{amsmath,amsthm, epsfig,times,graphicx,color}
\usepackage{amsfonts, graphicx,color,multirow}
\usepackage{amsfonts, graphicx,color,multirow}
\def\boxit#1{\vbox{\hrule\hbox{\vrule\kern6pt\vbox{\kern6pt#1\kern6pt}\kern6pt\vrule}\hrule}}

\newcommand{\PP}{\mathbb{P} }
\newcommand{\MM}{\mathbb{M} }
\newcommand{\GG}{\mathbb{G} }

\newcommand{\infinity}{{\infty}}

\newtheorem{theorem}{Theorem}[section]

\newtheorem{lemma}{Lemma}[section]

\newtheorem{prop}{Proposition}[section]

\newtheorem{defn}{Definition}[section]

\def\argmax{\mbox{argmax}}
\def\argmin{\mbox{argmin}}

\def\labda1{\lambda_1}
\def\labda2{\lambda_2}

\def\argmin{\mbox{argmin}}

\def\comment#1{\relax}

\def\=in{\mathop{\rm =}}

\newcommand{\ol}[1]{\overline{#1}}
\newcommand{\ep}[0]{\mathbb{P}}

\newcommand{\bbP}{\mathbb{P} }

\newcommand{\beq}{\begin{eqnarray}}
\newcommand{\eeq}{\end{eqnarray}}
\newcommand{\beqn}{\begin{eqnarray*}}
\newcommand{\eeqn}{\end{eqnarray*}}

\begin{document}

\def\adrui#1{\vskip 2mm\boxit{\vskip 2mm{\color{blue}\bf#1} {\color{red}\bf --By Rui  \vskip 2mm}}\vskip 2mm}

\begin{frontmatter}

\title{ Asymptotics for change-point models under varying degrees of mis-specification}

\runtitle{Change-point mis-specification}

 \begin{aug}
\author{\fnms{Rui} \snm{Song}\thanksref{t1}\ead[label=e1]{rsong@ncsu.edu}},
\author{\fnms{Moulinath}  \snm{Banerjee}\thanksref{t2}\ead[label=e2]{moulib@umich.edu}}
\and
\author{\fnms{Michael} \snm{R.} \snm{Kosorok}\thanksref{t2}\ead[label=e3]{kosorok@unc.edu}}


\thankstext{t1}{Supported in part by Grant NSF-DMS 1007698 and 1309465.}
\thankstext{t2}{Supported in part by Grant NSF-DMS 1007751 and 1308890.}
\thankstext{t3}{Supported in part by Grant NCI P01 CA142538.}
\thankstext{t4}{We thank the editor, the AE and three referees for providing input and suggestions that vastly improved the quality of the paper.}
\runauthor{R. Song, M. Banerjee and M. R. Kosorok}

\affiliation{North Carolina State University, University of Michigan and University of North Carolina}

\address{Rui Song\\
Department of Statistics,\\
 North Carolina State University,\\
Raleigh NC, U.S.A.\\
\printead{e1}}

\address{Moulinath Banerjee\\
Departments of Statistics \\
University of Michigan \\
Ann Arbor, MI, U.S.A.\\
\printead{e2}}

\address{Michael R. Kosorok\\
Departments of Biostatistics \\
~~and Statistics \& Operation Research,\\
  University of North Carolina at Chapel Hill,\\
Chapel Hill, NC, U.S.A.\\
\printead{e3}}

\end{aug}


\begin{abstract}

Change-point models are widely used by statisticians to model drastic changes in the pattern of 
observed data. Least squares/maximum likelihood based estimation of change-points leads to curious 
asymptotic phenomena. When the change--point model is correctly specified, such estimates 
generally converge at a fast rate ($n$) and are asymptotically described by minimizers of a jump process. Under 
complete mis-specification by a smooth curve, i.e. when a change--point model is fitted to data described by a smooth curve, 
the rate of convergence slows down to $n^{1/3}$ and the limit distribution changes to that of the minimizer of a continuous Gaussian 
process. In this paper we provide a bridge between these two extreme scenarios by studying the limit behavior of change--point 
estimates under varying degrees of model mis-specification by smooth curves, which can be viewed as local alternatives. We find that the 
limiting regime depends on how quickly the alternatives approach a change--point model. We unravel a family of `intermediate' limits that 
can transition, at least qualitatively, to the limits in the two extreme scenarios. The theoretical results are illustrated via a set of carefully designed 
simulations. We also demonstrate how inference for the change-point parameter can be performed in absence of knowledge of the underlying scenario by 
resorting to subsampling techniques that involve estimation of the convergence rate. 
\end{abstract}

\begin{keyword}[class=MSC]
\kwd[Primary ]{62G20}
\kwd{62G05}
\kwd[; secondary ]{62E20}
\end{keyword}

\begin{keyword}
\kwd{Change-point}
\kwd{Model mis-specification}
\end{keyword} 
\end{frontmatter}

\section{Introduction}
The study of change-point models has a long and rich history in the statistics and econometrics literature. Change-point models, where a signal function shows an abrupt transition at one or more points in its domain can be used to study phenomena that are subject to sudden shock effects, or which show natural phase-transitions at different stages of evolution. Applications are many and varied and arise in the analysis of climate data \citep{Lund02}, estimation of mixed layer depth from oceanic profile data \citep{Thomson03}, structural breaks in economics \citep{BP:98,BP:03}, quality control and dynamical systems in an engineering context \citep{Lai:95}, and genetics \citep{SZ:12}, to name a few.
Sequential methods for change-point detection have been around for a very long time; the literature here is truly huge, with a comprehensive treatment in the book by \cite{BN:93} and the excellent review paper by \cite{Lai:01}, but see also \cite{csorgo1997limit} which has an in-depth study of limit theorems in change--point analysis. On the other hand, inference on jump-discontinuities (change-points) in an otherwise smooth curve based on observed or designed data has also received attention in the nonparametric as well as the survival analysis literature: see, for example, \cite{GHK99,HM:03,Koso:Song:Inf:06}, \cite{LBM09,Loader96}, \cite{Muller92,MS:97} , \cite{Pons:esti:2003,Ritov90} and references therein. A canonical change-point model which illustrates many important features of this genre of problems is given by:
$$Y = \beta_l1(X\le \theta) + \beta_u1(X > \theta) + \epsilon, $$
where the predictor $X \in [0,1]$ is assumed to be a continuous random variable, $\beta_l \neq \beta_u$ are fixed constants, $\epsilon$ is
a continuous random variable, independent of $X$ with zero expectation and finite
variance. The parameters of interest are the change-point parameter $\theta$ and the regression parameter $(\beta_l, \beta_u)'$.
For this model, the least squares estimator of the change-point parameter converges to the truth at rate $n$, with the limit distribution being described by the minimizer of a two-sided, compound Poisson process.  The asymptotic distribution of the least squares estimates of $(\beta_l, \beta_u)'$ are normal, and are unaltered by estimation of the change-point: i.e., they have the same distribution as the least squares estimates that would have been obtained if $\theta$ were known. The detailed analysis can be found in \cite{Kta06}. A closely related model allows the parameters $\beta_l$ and $\beta_u$ to approach each other with increasing sample size $n$ (as opposed to staying fixed in the above display). As long as $\beta_u - \beta_l$ approaches 0 at a rate slower than $n^{-1/2}$, the change-point can be estimated. However, due to the loss of signal in this model, the rate of convergence of the LSE of $\theta$ slows to $n^{1-2\,\xi}$, where
$\beta_u - \beta_l = O(n^{-\xi})$; furthermore, the limit distribution is now starkly different and described by the minimizer of Brownian motion plus triangular drift. See for example, \cite{BB:76} for an early treatment of this problem, and \cite{MS:97} for a nonparametric incarnation. \cite{Huskova99} considered estimators in location models with various gradual changes and showed that the limit behavior of least--squares type estimators of the change point in these models depends on the type of gradual change. 

A natural question, as in most statistical problems, is the effect of mis-specification on the change-point estimator. Suppose first that the true model is of the form $Y = f(X) + \epsilon$ where $f$ is actually smooth but that the model $\beta_l\,1(X \leq \theta) + \beta_u\,1(X > \theta)$ is fitted instead. This is what happens, for example, in CART where the change-point analysis represents the best approximation of a binary decision tree (piecewise constant function with a single jump, also called a stump) to $f$. \cite{BY:02} and \cite{BM:07} studied the asymptotics of the estimates of the change-point and the regression coefficient in this problem and showed that in this setting cube-root asymptotics with Chernoff limit distributions obtain. As shown by \cite{BM:07}, all three least squares estimates converge at the slower $n^{1/3}$ rate because the change-point estimation depends on local features of the smooth regression curve, which are more complex in comparison to when the true regression function is a stump model. Therefore, change-point estimation and inference are highly unstable under model mis-specification by a smooth curve due to this sharp fall in the estimator's rate of convergence: from a rate as fast as $n$ under the true change-point model, to only a cube-root rate under a smooth curve.

While this is an interesting finding, this formulation does not quite capture the more subtle issue of how the degree of mis-specification affects the estimates of the stump parameters: to elaborate, consider functions $f_1$ and $f_2$ that are both smooth, but suppose that one is linear and the other a sigmoidal function with a sharp ascent. Clearly $f_2$ is much closer to a stump-model than $f_1$, so fitting the mis-specified change-point model should be less consequential in the case of $f_2$ than $f_1$. But the fixed-model approach described in the above paragraph does not satisfactorily capture this issue. This motivates us, in this work, to consider models where the degree of mis-specification is allowed to change -- diminish, in fact -- as $n \rightarrow \infty$ and to explore the consequences of this diminishing mis-specification on the behavior of the stump estimates. In particular, how does the rate of mis-specification bear upon this behavior in terms of rates of convergence and limit distributions? 

Our strategy considers a sequence of models $Y = f_n(X) + \epsilon$, where $f_n$ converges to a stump function at a rate controlled by a parameter $\alpha_n \rightarrow \infty$. We find that if the $f_n$'s converge to a stump slowly enough ($\alpha_n = o(n)$), the  limit distribution of the change-point estimator stays identical to the case $f_n \equiv f$, the fixed function setting of \cite{BM:07}, though the rate of convergence can be accelerated to (almost) $n$; if $f_n$'s approach the stump rapidly ($n = o(\alpha_n)$), the rate and limit distribution are identical to those that obtain when the true function is a stump, whereas, at the boundary $n \sim \alpha_n$, the limit distribution is different from either of the previous scenarios and belongs to a family of distributions that can transition, in a manner to be made precise in Section 4, to the Chernoff distribution (the limit with $\alpha_n = o(n)$) on one end and the minimizer of a compound Poisson process (the limit distribution with $n = o(\alpha_n)$) on the other. The joint limit behavior of the estimates of the levels of the stump and the jump-point, however, show an abrupt change as one changes from
$\alpha_n = o(n)$ to $n = O(\alpha_n)$: in the former case, the normalized estimates are asymptotically correlated with correlation 1 (i.e. linear functions of one another), which is also what happens in the ``fixed $f$'' scenario, while for the latter the estimates of the levels are asymptotically independent of that of the change-point. Viewing these $f_n$'s as a sequence of local alternatives to the limiting null model, a stump, the above phenomena are qualitatively identical to what transpires with the MLE in regular parametric models under a sequence of local alternatives, depending on how quickly the alternatives approach the null, an analogue we discuss more fully in the 
final section. 

The problem addressed in this paper should be contrasted with the `local alternative'--type models considered in \cite{BB:76}, \cite{MS:97} and \cite{Huskova99}. In all these papers, the limit of the sequence of change--point models considered-- the so-called `null' model -- is a smooth function \emph{without} a change--point, whereas we have the \emph{reverse} scenario: our sequence of models are smooth functions that, in the limit, produce a discontinuous change--point model. Our work is also quite different from inference in settings where the change-point is not a discontinuity but represents a point of smooth and/or gradual change; see, for example, \cite{VD15}. 
To the best of our knowledge, our work is the first attempt at providing a comprehensive as well as systematic understanding of the behavior of change-point models under local smooth alternatives. We hope that it will stimulate more investigation into this relatively uncharted territory. 

An interesting aspect of our approach is the fact that our paradigm approximates fixed function settings as the sample size varies. Consider again, the function $f_2$ introduced in the previous page, with its sharp ascent. As we show in our simulations (Section 6) , with a fixed steep function (not changing with $n$), each of our three asymptotic regimes takes turns being the best approximation to the sample distribution of the estimated changepoint. Specifically, for small sample sizes, the steep function is indistinguishable from the stump model and the compound Poisson process is the best approximation; for moderate sample sizes, the intermediate regime is the best approximation; while for large sample sizes, the fact that the function is not a stump is detectable by the estimator and the Chernoff limit becomes the best approximation. In other words, the proposed contiguous model sequence realistically approximates the range of distributional behavior which can be found in practical data settings.

From the perspective of statistical inference, the key contribution of our work is the formulation of a concrete theoretical framework for adaptive estimation of the change-point parameter under possible mis-specification through a subsampling procedure, elaborated in Section 5.  Consider the following basic inference problem: 
\emph{Given data $\{Y_i, X_i\}$ from a regression model to which a stump model has been fit, how would a practitioner go about constructing a confidence interval for the change-point parameter? The true  underlying (unknown) model: $Y = v_0(X) + \epsilon$ may not even be a change-point model but potentially in the proximity of one.} The framework of our paper then allows a way of making inference on the change-point parameter by (i) \emph{couching} the underlying model in a sequence of models $f_n$ determined by an unknown $\alpha_n$, (ii) providing a meaningful interpretation to the population change-point parameter, and (iii) last but not least, building confidence intervals for the parameter by adaptively estimating the correct regime for the given data-set through a subsampling procedure which teases out the correct degree of calibration (the convergence rate) from the data itself.

Before we proceed further, a word of clarification regarding the use of the term `change point' is in order. It is true that the underlying regression functions in the framework 
of this paper are smooth and therefore do not possess a change-point in the conventional sense of the term. The change-point model is used as a working model to fit the 
data; in that sense, the change-point is really a split-point in the spirit of \cite{BM:07}. However, if one takes the point of view that any change-point model fitted by a statistician is really an approximation to some underlying continuous model, with better fits corresponding to continuous models that are close to a model with a discontinuity -- so the notion of a change-point is really a convenient idealization -- then the term `change-point' can be used without further scope of confusion. Indeed, the title of the 
paper emphasizes this view in highlighting the mis-specification angle upfront! 

The rest of the paper is organized as follows. In Section 2, we formulate the regression problem and list our assumptions with interpretations. In Section 3, we systematically establish the asymptotic results including consistency, convergence rates and weak convergence of the change--point estimator. The connections among the different limiting distributions obtained for different choices of $\alpha_n$ are established in Section 4. Section 5 describes the adaptive inference strategy. Section 6 presents empirical evidence from a simulation study where a stump model is fit to data arising from a smooth regression curve and illustrates when the different asymptotic regimes comes into play: this is seen to depend on the nature of the underlying smooth function as well as the sample size, and the boundary case is seen to provide more robust approximations than the others. Section 7 discusses connections of our results to a number of other problems and scope for furthering this direction of research. 

{





\section{Change-point Models under Model Misspecifications}
\subsection{The Model Set-up}
At stage $n$, the observed data $(Y_i,X_i)$, $i=1,\ldots,n$ are $n$ i.i.d. copies of $(Y,X)$, where
$Y = f_n(X) + \epsilon$: here $E(\epsilon) = 0$, $\epsilon$ is independent of $X$ with bounded second moment
and $X$ follows some distribution $F_X$ on $[0,1]$. Thus, we have a sequence of models (changing with $n$).  The functions $f_n$ will be constructed to be smooth but will converge to a stump function as described later in this section. 

At each stage $n$, our working model will be a stump of the form $g(x;\theta, \beta_l, \beta_u) = \beta_l\,1(x \leq \theta) + \beta_u\,1(x > \theta)$ and the best working model will be determined from the sample via least squares. Denote a generic $(\theta, \beta_l, \beta_u)$ by $\psi$, and let
\begin{eqnarray*}
\hat \psi^n \equiv (\hat \theta^n, \hat \beta_l^n, \hat \beta_u^n)^T &= &  \mbox{argmin}_{\psi}
\sum_{i=1}^n\,(Y_i - g(X_i, \psi))^2 \\
& = &
\mbox{argmin}_{\psi}\,
\mathbb{P}_n [y - g(x; \psi)]^2\},
\end{eqnarray*}
with $\mathbb{P}_n$ denoting the empirical measure of the data $\{Y_i,X_i\}_{i=1}^n$. Letting $P_n$ denote the true distribution
at stage $n$, the corresponding population parameter $\psi^n=( \theta^n,\beta_l^n, \beta_u^n)^T$
is defined through the least squares estimation problem:
\begin{eqnarray*}
\psi^n &=& \argmin_{\psi} M_n(\psi) \\
&=&
\argmin_{\psi}P_n m_{\psi}(X,Y)=\argmin_{\psi}P_n\{Y - \beta_l 1(X \le \theta) - \beta_u 1(X >\theta) \}^2.
\end{eqnarray*}
We assume that there is a unique (population) minimizer $\psi^n \equiv (\theta^n, \beta_l^n, \beta_u^n)$, with $\beta_l^n \neq \beta_u^n$ at stage $n$.

To focus on the main ideas, we consider a sequence $f_n$ of the type:
$$f_n(x) = f(\alpha_n(x-\theta^0)),
$$
where $f$ is a smooth bounded monotone (increasing) function defined on $\mathbb{R}$ and $\alpha_n$ a sequence going to $\infty$. Denote $f(-\infty)$ by $\beta_l^0$ and $f(\infty)$ by $\beta_u^0$. As $n$ goes to infinity, $f_n(x)$ then converges to the stump
$$f_0(x) \equiv \beta_l^0 1(x \leq \theta^0) + \beta_u^0 1(x > \theta^0)$$
at all points except $\theta_0$. We let $\psi^0=(\theta^0, \beta_l^0, \beta_u^0)^T$ denote this limiting
population parameter. Note that the speed of convergence of $f_n$ to $f_0$ is regulated by the parameter $\alpha_n$.
Define $\xi^n = \alpha_n(\theta^n - \theta^0)$, which can be viewed as a rescaled ``bias'' term due to model misspecification.

Let $\xi^n = \alpha_n(\theta^n-\theta^0)$ and $\xi^0 = f^{-1}((\beta_l^0 + \beta_u^0)/2)$. From Theorem 2.2 below, $\lim_{n \rightarrow \infty} \xi^n = \xi^0$.

From the normal equations that characterize $\psi^n$, we have
\beqn
&&\xi^n=f^{-1}((\beta_l^n + \beta_u^n)/2),\\
&&\beta_l^n = \frac{\int_0^{\theta^n} f(\alpha_n(x-\theta^0))p_X(x) dx }{P(X\le \theta^n)},\\
&&\beta_u^n = \frac{\int_{\theta^n}^1 f(\alpha_n(x-\theta^0))p_X(x) dx }{P(X >  \theta^n)},
\eeqn


Although the working model is an oversimplification of the true model at each fixed $n$, as $n$ gets larger, the approximation to the true model is better. It will be shown later that $\psi^n$ converges to its limit $\psi^0$. The statistic $\hat{\psi}_n$ defined earlier estimates $\psi^n$ and therefore, indirectly, $\psi_0$. We note here that the minimizer of $\mathbb{P}_n [y - g_n(x; \psi)]^2$ is not necessarily unique, so in the case of multiple minimizers, we take $\hat{\psi}_n$ to be the minimizer with the smallest value of the first co-ordinate (if two minimizers have identical first co-ordinates their last two co-ordinates must also coincide). For simplicity of reference, we call this the \emph{smallest} argmin. We will study the asymptotic behavior of $\hat \psi^n$ as $\alpha_n$ converges to infinity at different rates.

\subsection{Assumptions}
We now describe our assumptions on the model formulated above. 


\begin{itemize}
\item [A.] $f(x)$ is continuously differentiable in an open neighborhood $N$ of $\xi^0$ with $f'(\xi_0)>0$. 

\item [B.] The density $p_X(x)$ does not vanish and is continuously differentiable in a neighborhood of $\theta_0$.



\item[C1.]
\beqn
\inf_n\liminf_{|h_1|\rightarrow \infty} \frac{1}{|h_1|}\int_0^{h_1}[f(t+\xi^n) - f(\xi^n)]dt>0.
\eeqn

\item[C2.] For any positive constant $K$, $\inf_n\inf_{u \in [-K+\xi^n, K+\xi^n]} f'(u)>0$.
\item [C3.] The integrals $\int_{-\infty}^{\xi^0}(f(t) - \beta_l^0)dt$
and $\int_{\xi^0}^{\infty}(f(t) - \beta_u^0)dt$ exist and are denoted as $\xi^l$ and $\xi^u$ respectively.
\end{itemize}

Assumptions A and B are adapted from the conditions in Banerjee and McKeague (2007). Assumption C1 says that the average increase of $f$ over all sufficiently large finite intervals with $\xi^n$ as an end-point should be bounded away from 0. Assumption C2 is, essentially, a positivity condition on the derivative of $f$ in every compact neighborhood of $\xi_0$. Assumption C3 figures in calculating the asymptotic bias of $\beta_l^n$ and $\beta_u^n$ for $\beta_l^0$ and $\beta_u^0$ respectively. Note that this assumption implies that $\int_{-\infty}^0\{f(t)-\beta_l^0 \}^2dt$ and $\int_0^{\infty}\{f(t)-\beta_u^0 \}^2dt$ are both $O(1)$.

In the sequel, it should be understood that the proof of any lemma, proposition or theorem that does not appear in the main text has been relegated to the appendix. 
\subsection{Limiting behavior of $\psi^n$}
We establish the asymptotic behavior of $\psi^n$, the stage $n$ population parameter in two steps. First, we show the consistency of $\psi^n$ for $\psi^0$ and next, we establish the convergence rates and calculate the limiting (normalized) bias of $\psi^n$ for $\psi^0$. Note that the convergence results of $\psi^n$ to $\psi^0$ in all three steps are deterministic.
The following theorem establishes the consistency of $\psi^n$.
\begin{theorem}\label{the-det-csist}
Under Assumptions A, B, and C, $\lim_{n \rightarrow \infty}(|\theta^n - \theta^0| + |\beta_l^n - \beta_l^0|+|\beta_u^n - \beta_u^0|)= 0$.
\end{theorem}
The next theorem deals with convergence rates and asymptotic bias.
\begin{theorem}\label{the-det-rate}
Under Assumptions A--C,  $\lim_n \alpha_n(\theta^n - \theta^0)=\xi^0$,
$\lim_n\alpha_n(\beta_l^n - \beta_l^0)=p_X(\theta^0)F_X(\theta^0)^{-1}\xi^l$, and
$\lim_n\alpha_n(\beta_u^n - \beta_u^0)=p_X(\theta^0)F_X(\theta^0)^{-1}\xi^u$.
\end{theorem}


{\bf Remark:} Since $\alpha_n(\theta^n - \theta_0) = f^{-1}((\beta_l^n + \beta_u^n)/2)$, and by Theorem \ref{the-det-csist}, $\beta_l^n$ and $\beta_u^n$ converge to $\beta_l^0$ and $\beta_u^0$ respectively, it is immediate that $\alpha_n(\theta^n - \theta_0)$ converges to $\xi_0$ as defined earlier in this section. The proofs of the remaining two convergences are given in the appendix. 

\section{Asymptotic Results for $\hat {\psi}_n$}
We will present the asymptotic results for $\hat {\psi}_n$ in three subsections in the order of consistency, convergence rates and weak convergence.
\subsection{Consistency}
We first establish Euclidean consistency for $\hat{\psi}^n$, where the results are summarized in Theorem \ref{the-con}.

\begin{theorem}\label{the-con}
Under Assumptions A--C, $|\hat \theta^n - \theta^n|+|\hat \beta_l^n-\beta_l^n|+|\hat \beta_u^n-\beta_u^n|=o_{P}(1)$.
\end{theorem}

\subsection{Rate of convergence}
In this section, we establish the convergence rates for change-point estimators under different degrees of model misspecification.
As an important first step, we introduce a dichotomous distance that describes the variation of the population criterion function about its minimizer.
\beqn
\rho_1(\psi, \psi^n) = \max\{ \alpha_n^{1/2}\,|\theta - \theta^n|, ~|\beta_l - \beta_l^n|, ~|\beta_u - \beta_u^n|\},
\eeqn
\beqn
\rho_2(\psi, \psi^n) = \max \{|\theta - \theta^n|^{1/2}\,,~ |\beta_l - \beta_l^n|, ~|\beta_u - \beta_u^n| \}.
\eeqn
The following lemma is about a unified distance which enables a certain expansion of the objective function.

\begin{lemma}\label{lem-dis}
Under Assumptions A--C, it follows that for $\psi$ in a neighborhood of $\psi^n$ defined as :$\{ \psi: d_n(\psi,\psi^n) < \delta_0\}$ for some small $\delta_0>0$, there exists a positive constant $E_0$ such that 
\beq\label{e-dicho}
M_n(\psi) - M_n(\psi^n) \geq E_0\,d_n^2(\psi, \psi^n),
\eeq
where ${d}_n^2(\psi,\psi^n)= \rho_1^2(\psi, \psi^n)1(|\theta - \theta^n| \leq K\,\alpha_n^{-1}) +
\rho_2^2(\psi, \psi^n)1(|\theta - \theta^n|
>  K\,\alpha_n^{-1}).$
\end{lemma}


This dichotomous nature of the distance $d_n$ is really what drives the convergence rate of $\hat \psi_n$. It reflects the fact that the magnitude of the fluctuation of $M_n$ around $\psi^n$ is governed by where $\theta$ falls with respect to a (shrinking) $\alpha_n^{-1}$ order neighborhood of $\theta^n$.
If $\theta$ falls in the shrinking neighborhood, the growth of $M_n$ around $\psi^n$ in the first co-ordinate is at least of order $\alpha_n(\theta - \theta^n)^2$; if not, the growth is at least of order $|\theta - \theta^n|$, which appears in the classic correctly specified change-point problem considered in Kosorok (2008).  Note that the order of $\alpha_n(\theta - \theta^n)^2$ is dominated by that of $|\theta - \theta^n|$, precisely when $|\theta - \theta^n|$ is $O(1/\alpha_n)$: this \emph{slower} growth of $M_n$ in its first co-ordinate in a shrinking Euclidean neighborhood is what makes the convergence rate fall below $n$ for a wide range of $\alpha_n$. For slow--growing $\alpha_n$, which can be considered as $\alpha_n$ essentially behaving like a constant, we converge towards the $\rho_1$ setting and the distance function of Banerjee and McKeague (2007) and approach the $n^{1/3}$ convergence rate for $\hat{\theta}^n$ obtained in their work; for rapidly growing $\alpha_n$, we move towards the $\rho_2$ setting and the distance function in Kosorok (2008), and, approach the $n$-rate of convergence instead. The precise statements of the convergence rates appear in Theorem 
\ref{the-rates} below. 
\newline
\newline
We next calculate bounds on the modulus of continuity of the empirical process with respect to this distance: this is one of the key ingredients that dictates the convergence rate. The dichotomous nature of the
distance requires exercising some care via calculating
\[ P_n^{\star}[\sup_{d_n(\psi,\psi^n)<\delta} |\GG_n\{m_{\psi}(X,Y) - m_{\psi^n}(X,Y)\}|],\]
where $\GG_n m(\cdot) = (\PP_n -P_n)m(\cdot)$ for function $m(\cdot)$.
By definition of the distance ${d}_n(\psi,\psi^n)$, for some $\delta>0$ we have 
\beqn
&& \{{d}_n(\psi,\psi^n) < \delta\} \\
&&= \{ \rho_1(\psi, \psi^n)< {\delta}, |\theta - \theta^n| \leq 1/\alpha_n\} \cup
\{\rho_2(\psi, \psi^n)   <  {\delta}, |\theta - \theta^n|^{1/2} >1/\sqrt{\alpha_n}\} .
\eeqn
For $\delta \leq 1/\sqrt{\alpha_n}$, the second term on the right side is the null set and since for this range $\delta/\sqrt{\alpha_n} \leq 1/\alpha_n$,
we have $\{{d}_n(\psi,\psi^n) < \delta\} = \{\rho_1(\psi, \psi^n) < \delta\}$.
\newline
On the other hand, for $\delta > 1/\sqrt{\alpha_n}$,
$\delta/\sqrt{\alpha_n} > 1/\alpha_n$ and the set $\{{d}_n(\psi,\psi^n) < \delta\} =
\{\rho_2(\psi,\psi^n) < \delta\}$. In the next lemma we establish the order of modulus of two function classes which will be used for the convergence rates, as stated in Theorem \ref{the-rates}.


\begin{lemma}\label{lem-mod}
Under Assumptions A--C, we have that for $0< \delta \leq 1/\sqrt{\alpha_n}$,
\beq\label{e-small}
E_n^{\star}[\sup_{d_n(\psi,\psi^n)<\delta} |\GG_n\,(m_{\psi}(X,Y) - m_{\psi^n}(X,Y))|\;\;] \lesssim \,\frac{\delta^{1/2}}
{\alpha_n^{1/4}} ,
\eeq
where $E_n^{\star}$ denote the outer expectation at stage $n$.
On the other hand, for $\delta > 1/\sqrt{\alpha_n}$,
\beq\label{e-large}
E_n^{\star}[\sup_{d_n(\psi,\psi^n)<\delta} |\GG_n\,(m_{\psi}(X,Y) - m_{\psi^n}(X,Y))|\;\;] \lesssim \delta .
\eeq
\end{lemma}
{\bf Remark:} The proof of Lemma \ref{lem-mod} involves reasonably standard arguments that use maximal inequalities to control the expected modulus of continuity of an empirical process via the magnitude of an envelope function and an entropy integral. The proof of Lemma \ref{lem-dis} needs more careful handling; in particular, it requires analyzing the fluctuation of $M_n$ about $\psi^n$ in terms of two components: the fluctuation about the first co-ordinate keeping the others fixed plus the fluctuation about the second and third co-ordinates keeping the first fixed. This is formalized in Lemma 0.1 in the appendix, the key preparatory result for the proof of Lemma \ref{lem-dis}. 

\begin{theorem}\label{the-rates}
Under Assumptions A--C, we have
\begin{enumerate}
\item[i.] When $\alpha_n = o(n)$, $n^{1/3}\alpha_n^{2/3} |\hat{\theta}^n - \theta^n| + n^{1/3}\alpha_n^{1/6}|\hat {\beta}_l^n - \beta_l^n|+ n^{1/3}\alpha_n^{1/6}|\hat {\beta}_u^n - \beta_u^n|=O_{P}(1)$.
\item[ii.]  When $\alpha_n = n$, $n\,|\hat{\theta}^n - \theta^n| + \sqrt{n}|\hat \beta_l^n - \beta_l^n|+ \sqrt{n}|\hat \beta_u^n - \beta_u^n|=O_{P}(1)$.
\item[iii.] When $n = o(\alpha_n)$, $n\,|\hat{\theta}^n - \theta^n|+ \sqrt{n}|\hat \beta_l^n - \beta_l^n|+ \sqrt{n}|\hat \beta_u^n - \beta_u^n|=O_{P}(1)$.
\end{enumerate}
\end{theorem}
{\it Proof of Theorem \ref{the-rates}.}

From Lemma \ref{lem-mod}, we have for $\delta \leq 1/\sqrt{\alpha_n}$,
\[
E_n^{\star}[\sup_{{d}_n(\psi,\psi^n)<\delta} |\GG_n\,(m(x,y,\psi) - m(x,y,\psi^n))|\;\;] \lesssim \frac{\delta^{1/2}}{\alpha_n^{1/4}} \,.
\]
On the other hand, for $\delta > 1/\sqrt{\alpha_n}$,
\[ E_n^{\star}[\sup_{{d}_n(\psi,\psi^n)<\delta} |\GG_n\,(m(x,y,\psi) - m(x,y,\psi^n))|\;\;] \lesssim \delta \,.\]
To apply Theorem 0.2 in the appendix, we are then led to a bounding function $\phi_n(\delta)$
for the modulus of continuity which is given by
\[ \phi_n(\delta) = \frac{\delta^{1/2}}{\alpha_n^{1/4}}\,1\left(\delta \leq \frac{1}{\sqrt{\alpha}_n} \right) + \delta\,1\left(\delta > \frac{1}{\sqrt{\alpha_n}}
\right) \,.\]
It is easily seen that $\phi_n(\delta)/\delta^{\alpha}$ is a decreasing function for $\alpha=1$.
Solving $r_n^2\,\phi_n(1/r_n) \leq \sqrt{n}$ yields
\begin{equation}
\label{rate-equation}
\frac{r_n^{3/2}}{\alpha_n^{1/4}}\,1(r_n \geq \sqrt{\alpha_n}) + r_n\,1(r_n < \sqrt{\alpha_n}) \leq \sqrt{n} \,.
\end{equation}

Next we analyze the rate from (\ref{rate-equation}) via isolating three cases for different choices of $\alpha_n$ one by one.

For the first case, considering $\alpha_n$ going to $\infty$ but no faster than $n$, i.e. $\alpha_n = o(n)$, we seek a solution with $r_n \geq \alpha_n^{1/2}$. To see this, suppose $r_n < \sqrt{\alpha_n}$. Then the solution is $r_n = \sqrt{n}$. Therefore, $\sqrt{\alpha_n}>\sqrt{n}$.  This is a contradiction, however, since by our condition, $\sqrt{n}$ is eventually
larger than $\sqrt{\alpha}_n$. This leads to:
$r_n = (n^{1/2}\,\alpha_n^{1/4})^{2/3} = n^{1/3}\,\alpha_n^{1/6}$.  We hence conclude that:
\beqn
&&n^{1/3}\alpha_n^{1/6}\left\{\rho_1(\hat \psi_n, \psi_n)\,1(\alpha_n |\theta^{n} - \theta^n| \leq 1)
 + \rho_2(\hat \psi_n, \psi_n)
  \,1(\alpha_n |\theta^{n} - \theta^n| > 1) \right\}
\eeqn
is $O_p(1)$. This implies that
\beqn
n^{1/3}\alpha_n^{1/6}|\hat {\beta}_l^n - \beta_l^n|=O_P(1),~~n^{1/3}\alpha_n^{1/6}|\hat {\beta}_u^n - \beta_u^n|=O_P(1),
\eeqn
and that:
\[ n^{1/3}\alpha_n^{2/3} |\hat{\theta}^n - \theta^n|\,1(\alpha_n |\theta^{n} - \theta^n| \leq 1) + n^{2/3}\,\alpha_n^{1/3}|\hat{\theta}^n - \theta^n|1(\alpha_n |\theta^{n} - \theta^n| > 1) \]
is $O_p(1)$. Since $\alpha_n = o(n)$ it is
strictly slower than both $n^{1/3}\,\alpha_n^{2/3}$ and $n^{2/3}\,\alpha_n^{1/3}$, showing that $\alpha_n\,|\hat{\theta}^n - \theta^n|$ is $o_P(1)$. This then forces $1(\alpha_n |\hat \theta^{n} - \theta^n| > 1)$ to go to 0 in probability. Since this is a zero-one valued random variable, it is easily argued that the second term in the above display must converge to 0 in probability. Given any subsequence, we can find a further subsequence along which the indicator converges almost surely to 0, and is therefore identically 0 in the long run, whence the second term also has to be identically 0. We thus conclude that $n^{1/3}\alpha_n^{2/3} |\hat{\theta}^n - \theta^n|$ is $O_P(1)$.

For the second case, we consider $\alpha_n = cn$, for some positive constant $c$.  We note that $\alpha_n=n$ is equivalent to $\alpha_n=cn$ for any $c\in(0,\infty)$ since the $c$ can be absorbed into the function $f$ without loss of generality. From now on we will use $\alpha_n=n$ everywhere else. Both rates $n^{1/3}\,\alpha_n^{2/3}$ and $n^{2/3}\,\alpha_n^{1/3}$ are equal to $n$ and we conclude that $n\,|\hat{\theta}^n - \theta^n|$ is $O_P(1)$, $\sqrt{n}|\hat \beta_l^n - \beta_l^n|=O_P(1)$ and $\sqrt{n}|\hat \beta_u^n - \beta_u^n|=O_P(1)$.

For the third case, we consider $n = o(\alpha_n)$. In this case, the second part in (\ref{rate-equation}) becomes relevant i.e. we seek a solution with $r_n < \sqrt{\alpha}_n$. The $r_n$ from the first part ---$n^{1/3}\,\alpha_n^{1/6}$ --- is inconsistent with the condition that $r_n \geq \sqrt{\alpha_n}$. and we are led to the solution $r_n = \sqrt{n}$ which is indeed consistent with the condition $r_n < \sqrt{\alpha_n}$. Conclude that:
\[ (n\,\alpha_n)^{1/2} \rho_1^2(\hat \psi^n, \psi^n) 1(\alpha_n|\hat{\theta}^n-\theta^n| \leq 1) + n\,\rho_2^2(\hat \psi^n, \psi^n)1(\alpha_n|\hat{\theta}^n-\theta^n| > 1)=O_P(1) \,.\]
Since $n\,\alpha_n$ is faster than $n^2$, it follows that $n\, |\hat{\theta}^n - \theta^n|$ is $O_P(1)$,
$\sqrt{n}|\hat \beta_l^n - \beta_l^n|=O_P(1)$ and $\sqrt{n}|\hat \beta_u^n - \beta_u^n|=O_P(1)$.
On the other hand, by the observation that the least squares estimate $\hat{\theta}^n$ is at least as far as $\theta^n$ from the $X_i$ closest to the latter and the fact that this $X_i$  converges to $\theta^n$ at rate $n$ (in fact, $n\,|X_i - \theta^n|$ converges to an exponential distribution), it follows that $n$ must be the non-trivial rate of convergence. $\Box$

\subsection{Asymptotic distributions}
Having established the rate of convergence, we now determine the asymptotic distribution. In the following, we discuss three different cases. The first result is the asymptotic distribution for $\alpha_n = o(n)$, which follows a rescaled Chernoff distribution. Recall that Chernoff's distribution is the unique minimizer of $W(t) + t^2$ over all real $t$, where $W(t)$ is two--sided Brownian motion starting from 0. 

\begin{theorem}\label{the-weak-smalln}
Let $q_n = n^{1/3}\alpha_n^{1/6}(\alpha_n^{1/2},1,1)^T$ and $F_X(\cdot)$ be the cumulative distribution function of $X$. Denote the pointwise product on Euclidean space as ``$\circ$''.
Under Assumptions A--C, when $\alpha_n=o(n)$,
\begin{eqnarray*}
&&q_n \circ \left(\hat \theta^n - \theta^n, \hat \beta_l^n - \beta_l^n, \hat \beta_u^n - \beta_u^n \right) \rightarrow_d (1, c_1, c_2)\argmax_h Q(h), ~\mbox{where} ~Q(h) ~\mbox{has }\\
&&\mbox{a rescaled Chernoff distribution:~} Q(h) = aW(h) - bh^2,
\end{eqnarray*}
$W(\cdot)$ is a standard two-sided Brownian motion process on the real line, $a^2= \sigma^2p_X(\theta^0)$,
\beqn
b = \frac{1}{2} f'(\xi^0)p_X(\theta^0) - \frac{1}{8}(\beta_u^0 - \beta_l^0)p_X(\theta^0)^2\left(
\frac{1}{F_X(\theta^0)} +\frac{1}{1-F_X(\theta^0)} \right),
\eeqn
\beqn
c_1 = \frac{p_X(\theta^0)(\beta_u^0 - \beta_l^0)}{2F_X(\theta^0)},~\mbox{and}~
c_2 = \frac{p_X(\theta^0)(\beta_u^0 - \beta_l^0)}{2(1-F_X(\theta^0))}.
\eeqn
\end{theorem}
{\bf Remark:} Note the similarity of the above results to that in Theorem 2.1 of Banerjee and McKeague (2007). The regime $\alpha_n = o(n)$ can be interpreted as the \emph{slow regime} which yields asymptotic behavior similar to the situation in that paper where the smooth function $f_n \equiv f$ 
and does not change with $n$. The form of the limits is similar to those obtained in Theorem 2.1 but note the difference in convergence rates. While in Banerjee and McKeague (2007) the rate of convergence of all three parameters is $n^{1/3}$, in our current situation we do get an acceleration above this rate: for the change-point parameter, the accelerated rate can (almost) go up to $n$ and for the level parameters it can (almost) go up to $\sqrt{n}$ as $\alpha_n$ gets close to order $n$, these limiting rates being the rates of convergence for a correctly specified change--point model. Also note that the asymptotic correlation between the least squares estimate of the stump levels and that of the change--point is 1, whereas, in the cases to follow, these will be seen to be asymptotically independent. 
\newline
\newline
The next result is the asymptotic distribution for $\alpha_n = n$. This is the most interesting scenario and yields a new limit distribution. To deduce the limit distribution of
$\hat h_n=(\hat h_{1n},\hat h_{2n},\hat h_{3n})^T$, where
$\hat h_{1n} = n(\hat \theta^n - \theta^n)$,
$\hat h_{2n} = \sqrt{n}(\hat \beta_l^n-\beta_l^n)$ and
$\hat h_{3n} = \sqrt{n}(\hat \beta_u^n-\beta_u^n)$, we consider the limit of the process $h \mapsto Q_n(h)= n\bbP_n(m_{\psi_{n,h}} - m_{\psi^n})$,  where 
$$\psi_{n,h} = \psi^n + (h_1/n, h_2/\sqrt{n}, h_3/\sqrt{n})\;\;\;\mbox{and}\;\;\;  h=(h_1, h_2, h_3)^T\,.$$

The general scheme of argument runs as follows: We first derive a tractable approximation of $Q_n$, denoted $\tilde Q_n$, that is uniformly close to $Q_n$ in a sense to be made precise later. The advantage of $\tilde{Q}_n$ is its decomposability into three parts where each represents the contribution of a parameter. Next, the tightness of $\tilde Q_n$ is established, which coupled with finite--dimensional convergence furnishes the weak limit of $\tilde Q_n$. This, by the uniform closeness alluded to above is also the weak limit of $Q_n$. The final step involves deriving the weak convergence of the normalized estimators by the application of an appropriate continuous mapping theorem for the argmax/argmin functional. 
\newline
We start with the first step. From the results on convergence rates, we know that $\hat h_{n}=(\hat h_{1n}, \hat h_{2n}, \hat h_{3n})^T$ is uniformly tight and is the smallest argmin of $h \mapsto Q_n(h)= n\bbP_n(m_{\psi_{n,h}} - m_{\psi^n})$.  Observe that $m_{\psi_{n,h}}(X,Y) - m_{\psi_n}(X,Y)$
\beqn
&=&
2(Y - f_n(\theta^n))(\beta_u^n - \beta_l^n)\{1(X \le \theta^n+h_1/n) - 1 (X \le \theta^n ) \}\\
&&~~+ (2Y - 2\beta_l^n - h_2/\sqrt{n})1(X \le \theta^n + h_1/n)h_2/\sqrt{n}\\
&&+ (2Y - 2\beta_u^n - h_3/\sqrt{n})1(X > \theta^n + h_1/n)h_3/\sqrt{n}.
\eeqn

Consequently,
\begin{eqnarray*}
{Q}_n(h)&=&
2 (\beta_u^n - \beta_l^n)n \ep_n(Y - f_n(\theta^n))\{1(X \le \theta^n+h_1/n) - 1 (X \le \theta^n ) \}\\
&&+2\sqrt{n}\left[\ep_n(Y - \beta_l^n)1\{X\leq\theta^n+h_1/n\}\right]h_2 - \ep_n1\{X\leq\theta^n +h_1/n\}h_2^2\\
&&+2\sqrt{n}\left[\ep_n(Y - \beta_u^n)1\{X>\theta^n+h_1/n\}\right]h_3 - \ep_n1\{X>\theta^n+h_1/n\}h_3^2\\
&=&{T}_{1n}(h_1) + \hat{T}_{2n}(h_1,h_2)+\hat{T}_{3n}(h_1,h_3),
\end{eqnarray*}
where
\beqn
{T}_{1n}(h_1)=2 (\beta_u^n - \beta_l^n)n\ep_n(Y - f_n(\theta^n))\{1(X \le \theta^n+h_1/n) - 1 (X \le \theta^n ) \},
\eeqn
\beqn
\hat{T}_{2n}(h_1,h_2) = 2\sqrt{n}\left[\ep_n(Y - \beta_l^n)1(X\leq\theta^n+h_1/n)\right]u_2 - \ep_n1\{X\leq\theta^n +h_1/n\}h_2^2, ~\mbox{and}
\eeqn
\beqn
\hat{T}_{3n}(h_1,h_3)=2\sqrt{n}\left[\ep_n(Y - \beta_u^n)1(X>\theta^n+h_1/n)\right]u_3 - \ep_n1\{X>\theta^n+h_1/n\}h_3^2.
\eeqn

We now define $\tilde Q_n(h)$ as follows:
\beqn
\tilde{Q}_n(h)&&=T_{1n}(h_1) + 2\sqrt{n}\ep_n\left[\epsilon 1\{X\leq\theta^n\}\right]h_2 - \bbP_n\{X\leq\theta^n\}h_2^2\\
&&+2\sqrt{n}\ep_n\left[\epsilon 1\{X>\theta^n\}\right]h_3 - \bbP_n\{X>\theta^n\}h_3^2 \\
&&= T_{1n}(h_1) + T_{2n}(h_2)+T_{3n}(h_3) = \bbP_nT_1(h_1) + \bbP_nT_2(h_2) + \bbP_nT_3(h_3),~\mbox{where~}\\
&& T_{1}(h_1) =2n(\beta_u^n - \beta_l^n)(Y - f_n(\theta^n))\{1(X \le \theta^n+h_1/n) - 1 (X \le \theta^n ) \},\\
&& T_{2}(h_2) = 2\sqrt{n}\left[\epsilon 1\{X\leq\theta^n\}\right]h_2 - 1\{X\leq\theta^n\}h_2^2,~\mbox{and~}\\
&&T_{3}(h_3) = 2\sqrt{n} \left[\epsilon 1\{X>\theta^n\}\right]h_3 -  1\{X>\theta^n\}h_3^2.
\eeqn

In Lemma \ref{QQ} below, we show that $Q_n(h)$ and $\tilde Q_n(h)$ are uniformly close, as random elements in the space $\mathcal{D}_K$, where $\mathcal{D}_K$, $K \subset \mathbb{R}^3$ is the space of functions $q:~ K \mapsto \mathbb{R}$, $K$ being a compact rectangle in $\mathbb{R}^3$. Such functions $w(h_1, h_2, h_3)$ are piece-wise constant, hence, cadlag in the first argument, $h_1$, and are continuous in the last two arguments $(h_2, h_3)$. For each compact interval $C$ in $\mathbb{R}$, define $\Lambda_C$ to be the collection of continuous, strictly increasing
maps $\lambda: C \mapsto C$ such that $\lambda(C)=C$. Similar to \cite{SS:11}, define a norm on $\Lambda_C$ as follows: 
\begin{eqnarray}
\lambda \mapsto \|\lambda \| \equiv \sup_{s \neq t, s,t\in C}\left|\log \frac{\lambda(t)-\lambda(s)}{t-s}\right|.
\end{eqnarray}
Note that $K = I \times A$, necessarily, for a two-dimensional compact rectangle $A$ and a compact interval $I$. For $w_1$, $w_2\in \mathcal{D}_K$, we define the Skorohod topology as the one induced by the metric
\begin{eqnarray*}
d_K(w_1, w_2) \equiv \inf_{\lambda \in \Lambda_I}\left \{\sup_{u \in K}\left|
w_1(u_1, u_2, u_3) - w_2(\lambda(u_1), u_2, u_3) \right| + \|\lambda \| \right\}.
\end{eqnarray*}
Endowed with this metric, $\mathcal{D}_K$ is a complete separable metric space. 

\begin{lemma}\label{QQ}
Under conditions B--C2,
$Q_n - \tilde Q_n = o_{P}^K(1)$ in $(\mathcal{D}_K, d_K)$ for each $K$ above. The superscript $K$ in $o_{P}^K(1)$ indicates that the norm of the error is in terms of $d_K$.
\end{lemma}



To obtain the limit distribution of $Q_n(h)$, we next establish the uniform tightness of  $\{\tilde Q_n\}_{n=1}^{\infty}$. 
\begin{lemma}\label{lem-tight}
The process $\{\tilde Q_n\}_{n=1}^{\infty}$ is uniformly tight. 
\end{lemma}
We now define the limit process. Let $\{\nu^{+}(h): h \geq 0\}$ be a homogeneous Poisson process on $[0,\infty)$ with right continuous and left limit (in short RCLL) sample paths and rate parameter $p_X(\theta^0)$.
Let $\{\epsilon_i\}_{i=1}^{\infty}$ be i.i.d. versions of $\epsilon$ and distributed independently of $\nu^{+}(h)$. Let $S_i$ denote the time to the $i$'th arrival for the Poisson process $\nu^{+}$, i.e. $S_i = R_1
+ R_2 + \ldots + R_i$, where $\{R_j\}_{j=1}^{\infty}$ are the i.i.d. exponential inter-arrival times corresponding to $\nu^{+}(h)$. For $h \geq 0$, define:

\[\Lambda_1(h) =
\sum_{j=0}^{\nu^{+}(h)}\,\left(\epsilon_j + f(S_j + \xi^0) - f(\xi^0) \right) \,.\]

To define the process for $h \leq 0$, generate $\nu^{-}(h)$, an LCRR (left continuous with right limit) homogeneous Poisson process on $[0,\infty)$ with parameter $p_X(\theta^0)$ and $\{\tilde{\epsilon}_i\}_{i=1}^{\infty}$ i.i.d. $\epsilon$ again,
and independent of $\nu^{-}(h)$.  Also, $\nu^{-}$ and the $\tilde{\epsilon}_i$'s are generated independently of $\nu^{+}$ and $\epsilon_i$'s. Let $\tilde{S}_i$ denote the time to the $i$'th arrival for the
process $\nu^{-}$. For $h \leq 0$, define:
\[ \Lambda_1(h) =
\sum_{j=0}^{\nu^{-}(-h)}\,\left(-\tilde{\epsilon}_j + f(\xi^0) - f(-\tilde{S}_j + \xi^0)  \right) \,.\]
It can be easily seen that the process $\Lambda_1(h)$, thus defined, has independent increments.

We now show that on every compact rectangle $K$, $\tilde Q_n(u)$ converges to the tight process $Q(u)$ in the $d_K$ metric, where 
\begin{eqnarray*}
Q(u)&=& 2(\beta_u^0-\beta_l^0)\Lambda_1(u_1) + 2Z_1u_2+u_2^2P\{X\leq\theta^0\}+2Z_2u_3+u_3^2P\{X>\theta^0\}\\
&\equiv&2(\beta_u^0-\beta_l^0) \Lambda_1(u_1) + \Lambda_2(u_2)+ \Lambda_3(u_3),
\end{eqnarray*}
where $Z_1$ and $Z_2$ are mean zero independent Gaussians with respective variances $\sigma^2P\{X\leq\theta^0\}$ and $\sigma^2P\{X>\theta^0\}$ and $Z_1$, $Z_2$, and $\Lambda_1$ are all independent. The result is summarized in Theorem \ref{the-weak-n}.


\begin{theorem}\label{the-weak-n}
Under Assumptions A--C, when $\alpha_n = n$, the process $\tilde Q_n$ converges weakly to $Q$ in $\mathcal{D}_K$ for every compact rectangle $K$ 
in $\mathbb{R}^3$. Furthermore, via a continuous mapping argument, $\hat h_n \rightarrow_d h^{\star}$, where $$h^{\star} =(h_1^{\star}, h_2^{\star}, h_3^{\star})^T= \argmin_{h \in \mathbb{R}^3} Q(h)\,.$$ 
Also, $n(\hat \theta^n - \theta^n) = \argmin_{h_1}T_{1n}(h_1) + o_{P}(1)$ and converges weakly to $\hat \nu_{\Lambda_1}$, where $\hat \nu_{\Lambda_1} = \inf\{\nu:~\Lambda_1(\nu) = \min_{\nu} \Lambda_1\}$, while  $\sqrt{n}(\hat \beta_l^n - \beta_l^n)$ and $\sqrt{n}(\hat \beta_u^n - \beta_u^n)$ converge weakly to mean zero Gaussian variables with variances $\sigma^2/P(X \le \theta^0)$ and $\sigma^2/P(X > \theta^0)$ respectively. Finally, $n(\hat \theta_n-\theta^n)$, $\sqrt{n}(\hat \beta_l^n - \beta_l^n)$ and $\sqrt{n}(\hat \beta_u^n - \beta_u^n)$ are asymptotically independent.
\end{theorem}
{\bf Remark:} Note that, by the argmin of $Q$, we mean the \emph{smallest} argmin as with $M_n$ in Section 2, since there may be multiple minimizers 
with differing values of the first co-ordinate. 

The next result is the asymptotic distribution for $n = o(\alpha_n)$, when the rate of the rescaling parameter $\alpha_n$ going to infinity, i.e., the speed that the working model approaches the true model, is even faster than $n$. In this scenario, the obtained limiting distribution is identical with that obtained under correct specification: i.e. when the true model is $f_0(x; \psi)=\beta_l^0 1(x \le \theta^0) + \beta_u^0 1(x> \theta^0)$, the limit of the regression functions considered in this paper. The arguments for this case follow exactly the same pattern as the case $n = \alpha_n$, so we omit the details and only describe the limit process and the asymptotic convergence results. Note that the rate of convergence in the two cases: $\alpha_n = n$ and $n = o(\alpha_n)$ are identical, and $\hat{h}_n$ and $Q_n$ are therefore defined in the exact same way as for the case $\alpha_n = n$. 

Recall that $\{\nu^{+}(h): h \geq 0\}$ is a homogeneous Poisson process on $[0,\infty)$ with right continuous and left limit (RCLL) sample paths and rate parameter $p_X(\theta^0)$ and $\{\epsilon_i\}_{i=1}^{\infty}$ are i.i.d. versions of $\epsilon$ and distributed independently of $\nu^{+}(h)$. For $h \geq 0$, define:

\[\Lambda(h) =
\sum_{j=0}^{\nu^{+}(h)}\,\left(\epsilon_j + \beta_u^0 - f(\xi^0) \right) \,.\]

To define the process for $h \leq 0$, again consider $\nu^{-}(h)$ and $\{\tilde{\epsilon}_i\}_{i=1}^{\infty}$, exactly as defined before and independent of 
$\nu^+(h)$ and $\{\epsilon_i\}$. For $h \leq 0$, define:
\[ \Lambda(h) =
\sum_{j=0}^{\nu^{-}(-h)}\,\left(-\tilde{\epsilon}_j + f(\xi^0) - \beta_l^0  \right) \,.\]
It is easily seen that the process $\Lambda(h)$ has independent increments. Also, note that the process only depends on $f$ through its limits at $-\infty$ and $\infty$: this follows by recalling that $f(\xi_0) = (\beta_l^0 + \beta_u^0)/2$. The proof of the below theorem is skipped owing to its similarities to the 
proof of Theorem \ref{the-weak-n}. 
\begin{theorem}\label{the-weak-fastn}
Under Assumptions A--C, when $n=o(\alpha_n)$, $n(\hat \theta_n-\theta^n)$, $\sqrt{n}(\hat \beta_l^n - \beta_l^n)$ and $\sqrt{n}(\hat \beta_u^n - \beta_u^n)$ are asymptotically independent. Furthermore, $n(\hat \theta_n - \theta^n) = \argmin_{h}T_{1n}(h_1) + o_{P}(1)$ and converges weakly to $2(\beta_u^0 - \beta_l^0)\hat \nu_{\Lambda}$, where $\hat \nu_{\Lambda} = \inf\{\nu:~\Lambda(\nu) = \argmin \Lambda\}$, while $\sqrt{n}(\hat \beta_l^n - \beta_l^n)$ and $\sqrt{n}(\hat \beta_u^n - \beta_u^n)$ converge weakly to mean zero Gaussian variables with variances $\sigma^2/F_X(\theta^0)$ and $\sigma^2/(1-F_X(\theta^0))$ respectively.
\end{theorem}

\section{Connections among the different limit distributions}
The goal in this section is to explore the connections between the three limiting regimes that arise when considering the behavior of
$\hat{\theta}_n - \theta^n$ (appropriately normalized) for different values of $\alpha_n$. For $\alpha_n = o(n)$, we get Chernoff's
distribution, up to a constant, whereas minimizers of two-sided compound Poisson processes appear in the other two cases. For $\alpha_n = n$, the
limit distribution depends on the entire function $f$, whereas for $n = o(\alpha_n)$, the distribution depends only on the limiting change-point model $f_0$. We show below that the distribution in the intermediate case, $\alpha_n = n$, belongs to a family of ``\emph{boundary distributions}'' that can transition, at least qualitatively,  to each of the other two limits. For easy exposition, we first restrict attention to the following \emph{one-parameter}  version of our problem.  The case where $\beta_l$ and $\beta_u$ are unknown will be discussed later.

At stage $n$, consider the model $Y = f(\alpha_n\,(X - \theta_0)) + \epsilon$ with the levels $\beta_l^0$ and $\beta_u^0$ assumed \emph{known}.
 We estimate $\theta_0$ by 
 \[\ol{\theta}_n := \mbox{argmin}_{\theta}\,\PP_n[(Y - \beta_l^0)^2\,1(X \leq \theta) + (Y - \beta_u^0)^2\,1(X > \theta)] \equiv \mbox{arg} \min \,\MM_n(\theta)\,,\] 
 where $\MM_n(\theta) = \PP_n\,[(Y - 1/2)1(X \leq \theta)]$, the equivalence of the two criterion functions being a consequence of some simple algebra. As before, the smallest argmin is used.
 \newline
The population version of $\MM_n(\theta)$ is given by: $M_n(\theta) = P_n[(Y - 1/2)1(X \leq \theta)]$ and $\tilde{\theta}^n = \mbox{arg} \min_{\theta}\,M_n(\theta)$. As in the 3 parameter problem, let $\xi_0 = f^{-1}((\beta_l^0 +\beta_u^0)/2)$, let $a_0 = \sqrt{\sigma^2\,p_X(\theta_0)}$ and $b_0 = p_X(\theta_0)f'(\xi_0)/2$. It is not difficult to check that $\tilde{\theta}^n = \theta_0 + (1/\alpha_n)\,\xi_0$. The following theorem gives the distribution of $\ol{\theta}_n$ under the different regimes.

\begin{theorem}
\label{one-param-asym}
In the above one parameter model,
\begin{itemize}
\item[(a)] when $\alpha_n = o(n)$,
\[ n^{1/3}\,\alpha_n^{2/3}\,(\overline{\theta}_n - \tilde{\theta}^n) \rightarrow_d L \equiv \arg \min_h\;(a_0\,W(h) + b_0 \,h^2) \,;\]
\item[(b)] when $\alpha_n = n$,
\[ n\,(\overline{\theta}_n - \tilde{\theta}^n) \rightarrow_d \, \arg \min_h \,\Lambda(h)\,,\]
where \beqn
&& \Lambda(h) =  \left\{\sum_{j=0}^{\nu^{+}(h)}\,\left(\epsilon_j + f(\xi_0 + S_j) - f(\xi_0) \right)\right\}\,1(h \geq 0) \\
&&+
\left\{\sum_{j=0}^{\nu^{-}(-h)}\,\left(-\tilde{\epsilon}_j + f(\xi_0) - f(\xi_0 - \tilde{S}_j) \right) \right\}\,1(h < 0) \,,
\eeqn
where $S_j$'s and $\tilde{S}_j$'s are as defined previously;
\item[(c)] when $n = o(\alpha_n)$,
\[ n\,(\overline{\theta}_n - \tilde{\theta}_n) \rightarrow_d \, \arg \min_h \,\tilde{\Lambda}(h)\,,\]
where
\beqn
&& \tilde{\Lambda}(h) =  \left\{\sum_{j=0}^{\nu^{+}(h)}\,\left(\epsilon_j + \beta_u^0 - f(\xi_0) \right)\right\}\,1(h \geq 0) \\
&&+
\left\{\sum_{j=0}^{\nu^{-}(-h)}\,\left(-\tilde{\epsilon}_j + f(\xi_0) - \beta_l^0) \right) \right\}\,1(h < 0) \,.
\eeqn
\end{itemize}
\end{theorem}
{\bf Remark:} Note that the limit distributions in (b) and (c) are identical to those obtained for $n(\hat{\theta}_n - \theta^n)$ in the 3 parameter problem, while the limit distribution in case (a) is \emph{different}: the constant $b_0$ in the drift term is \emph{larger} than $b$ that shows up in the three parameter problem; see Theorem \ref{the-weak-smalln}. The smaller $b$ leads to a larger variance in the 3 parameter problem, the price of having to estimate the levels $\beta_l^0$ and $\beta_u^0$. In the settings (b) and (c), the estimation of the levels has no effect on the distribution of the change-point since the level estimates are asymptotically independent of the change--point estimate and therefore, the distributions in the 1--parameter and 3--parameter problems coincide. The proof of the above theorem is skipped as it involves easier versions of the arguments required to prove the distributional results in the 3 parameter problem. 

We now introduce a family of processes $\{\Lambda_c\}_{c>0}$ that generalizes the process $\Lambda$ appearing in the central case, (b). For $c > 0$, define:
\beqn
&& \Lambda_c(h) =  \left\{\sum_{j=0}^{\nu^{+}(h)}\,\left(\epsilon_j + f(\xi_0 + c\,S_j) - f(\xi_0) \right)\right\}\,1(h \geq 0) \\
&&+
\left\{\sum_{j=0}^{\nu^{-}(-h)}\,\left(-\tilde{\epsilon}_j + f(\xi_0) - f(\xi_0 - c\,\tilde{S}_j) \right) \right\}\,1(h < 0).
\eeqn
The parameter $c$ that dictates the above family is a scale parameter that regulates the shift of the increments of the generalized compound Poisson process $\Lambda_c$.
An instructive (statistical) way of thinking about $\Lambda_c$ is to consider the model: $Y = f_c(n(X-\theta_0)) + \epsilon$, with $f_c(t) \equiv f(ct)$. By calculations similar to those needed to prove Theorem \ref{one-param-asym}, we can show that:
\begin{equation}
\label{conv-c-central} 
n\,(\overline{\theta}_{n,c} - \tilde{\theta}^{n,c}) \rightarrow_d \, \arg \min_h \,\Lambda_c(h) \,,
\end{equation}
where $\overline{\theta}_{n,c}$ and $\tilde{\theta}^{n,c}$ are the analogues of $\overline{\theta}_n$ and $\tilde{\theta}^n$ in the one parameter model above, which corresponds to $c=1$.

The following results show that the distribution of the minimizer of $\Lambda_c$ approximates the limit distributions in the cases (a) and (c), as $c$ approaches 0 and $\infty$ respectively, for the one--parameter problem.
\begin{theorem}
\label{c-goes-to-zero}
Under Assumptions A--C, as $c \rightarrow 0$,
\[ c^{2/3}\,\arg \min_h\,\Lambda_c(h) \rightarrow_d L \equiv \arg \min_h\,[a_0\,W(h) + b_0\,h^2] \,.\]
\end{theorem}

\begin{theorem}
\label{c-goes-to-infinity}
Under Assumptions A--C, as $c \rightarrow \infty$,
\[ \arg \min_h\,\Lambda_c(h) \rightarrow_d \tilde{L} \equiv \arg \min_h \tilde{\Lambda}(h)\,.\]
\end{theorem}

Heuristically, Theorem \ref{c-goes-to-infinity} is somewhat easier to visualize. As $c \rightarrow \infty$, for every $j$, $f(\xi_0 - c\,\tilde{S}_j)$ goes to 0 almost surely and$f(\xi_0 + c\,S_j)$ to 1 almost surely, and by putting in these limiting values in the expression for $\Lambda_c$ we recover the process
$\tilde{\Lambda}$. This is not a rigorous verification, as we need to show that the convergence of the processes happens in a strong enough topology for distributional convergence of the argmin functional. This is accomplished in the proof of Theorem \ref{c-goes-to-infinity}. 
As far as Theorem \ref{c-goes-to-zero} is concerned, the crux of the argument lies in showing that an appropriately scaled version of $\Lambda_c$ (where scaling appears in the magnitude of the process as well as its argument) converges to a Brownian motion plus a quadratic drift; see Theorem 0.4 in \cite{Supple}. 

Define the sequence $c_n: = \alpha_n/n$. Consider first, case (c): $n = o(\alpha_n)$, where the statistical model can be written as $Y = f_{c_n}(n(X-\theta_0)) + \epsilon$ with $c_n \rightarrow \infty$. By (\ref{conv-c-central}), conclude that the distribution of $n(\overline{\theta}_n - \tilde{\theta}^n)$ can be approximated by that of $\mbox{arg} \min\,\Lambda_{c_n}(h)$. This, of course, is consistent with what we learn in Theorems \ref{c-goes-to-infinity} and \ref{one-param-asym}: as $c_n$ grows large in this case, by Theorem \ref{c-goes-to-infinity}, $\mbox{arg} \min\,\Lambda_{c_n}(h)$ and $\argmin_h\,\tilde{\Lambda}(h)$ are close in a distributional sense, and the latter is indeed the limit of $n(\overline{\theta}_n - \tilde{\theta}^n)$ in Case (c) of Theorem \ref{one-param-asym}. 

Next, consider case (a): $\alpha_n = o(n)$. As above, using (\ref{conv-c-central}), conclude that the distribution of $n\,(\overline{\theta}_n - \tilde{\theta}^n)$ can be approximated by that of $\arg \min\,\Lambda_{c_n}(h)$, as well. Since $c_n$ becomes small in this case, by Theorem \ref{c-goes-to-zero}, this can be approximated by $c_n^{-2/3}\,L$, which is essentially what Part (a) of Theorem \ref{one-param-asym} tells us. Thus, the family $\{\Lambda_c\}$ provides a uniform approximation to the limit distributions across the three different situations.

In the 3 parameter problem, when $\alpha_n = o(n)$, we know from Theorem \ref{the-weak-smalln} that
$n^{1/3}\,\alpha_n^{2/3}\,(\hat{\theta}_n - \theta^n) \rightarrow_d \arg \min_h\,(a_0\,W(h) + b\,h^2) := L'$, and $L'$ and $L$ have \emph{different}  distributions. The distribution of $n(\hat{\theta}_n - \theta^n)$ can then be approximated by that of $c_n^{-2/3}\,L'$. Noting that $L \equiv_d  (a_0/b_0)^{2/3}\,\mathbb{C}$ and $L' \equiv_d  (a_0/b)^{2/3}\,\mathbb{C}$, where $\mathbb{C} = \arg \min_h\,(W(h) + h^2)$ is the Chernoff random variable,  the distribution of $n(\hat{\theta}_n - \theta^n)$ can be approximated by that of $c_n^{-2/3}(b_0/b)^{2/3}\,L$, and therefore by $(b_0/b)^{2/3}\,\arg \min\,\Lambda_{c_n}(h)$. With $n = o(\alpha_n)$, it is not difficult to see that the distribution of $n(\hat{\theta}_n - \theta^n)$ in the 3 parameter case can still be approximated by  $\mbox{arg} \min_h\,\Lambda_{c_n}(h)$, as in the 1 parameter case.




\section{Adaptive inference for the change--point parameter}
Inference on $\theta^n$ when $\alpha_n$ is known can be achieved through subsampling or the ``m out of n'' bootstrap. To perform adaptive inference when $\alpha_n$ is \emph{unknown}, which is the case in practice, it is important to estimate it reliably. To this end, we resort to the results in \cite{subsampling99} who proposed a subsampling procedure when the convergence rate is unknown: the key idea is to use the data to first construct an estimate of the rate of convergence and then use this \emph{estimated rate} to produce a confidence interval for the parameter of interest. 
Following their idea, we describe an adaptive inference procedure for $\theta^n$ when $\alpha_n$ is \emph{unknown}. 

Consider $\alpha_n = n^{\gamma}$ where $0 \le \gamma < \infty$. (We restrict ourselves to this polynomial class as this covers essentially all interesting regimes and is 
tractable to deal with using the suggested method.)  By the asymptotic results of the previous section, we know that
$n^{1/3}n^{2\zeta/3}(\hat \theta_n - \theta^n) $ converges to a tight random variable, say $L$, where $\zeta = \gamma \wedge 1$, 
the minimum of $\gamma$ and $1$.  To construct a level $1-\alpha$ C.I. we proceed thus: 
\begin{itemize}
\item[(1)] Pick subsample sizes $n_1 < n_2 < n$ where $n_j = n^{\beta_j}$, with $1 > \beta_2 > \beta_1 > 0$. For each $j = 1,2$, collect a subsample of size $n_j$ 
without replacement $r$ times and run the change-point estimation procedure these $r$ subsamples to 
obtain change-point estimates $\{\hat \theta_{jk}^*\}_{k=1}^r$. 

\item[(2)] Next, note that for each $j$, the empirical distribution of the 
$\{n_j^{(1/3 +2\zeta/3)}(\hat \theta_{jk}^* - \hat \theta_n)\}_{k=1}^r$ (conditional on the given data) approximates the distribution of $L$. Using a moment approximation, we can then write: 
\[ \frac{1}{r}\,\sum_{k=1}^r\,n_j^{(1/3 +2\zeta/3)}|\hat \theta_{jk}^* - \hat \theta_n| \approx E(|L|)\;,\;\;j = 1,2 \,,\] 
and therefore: 
\[ \log E(|L|) \approx \log\,\left[ \frac{1}{r}\,\sum_{k=1}^r\,|\hat \theta_{jk}^* - \hat \theta_n|\right] + \frac{1+2\,\zeta}{3}\,\log n_j \;,\;\;j = 1,2.\]

\item[(3)] Equating the right-side of the above display for $j =1$ to that for $j=2$, a natural estimate of $\zeta$ is found by solving the equation:  
\beqn
\frac{1 + 2\hat \zeta}{3}=\left\{ \log\Bigl( \frac{n_2}{n_1}\Bigr)\right\}^{-1}\left[\log\left\{\frac{1}{r}\sum_k |\hat \theta_{1k}^*-\hat \theta_n| \right\}  - \log\left\{\frac{1}{r}\sum_k |\hat \theta_{2k}^*-\hat \theta_n| \right\} \right]\,.
\eeqn
This formula is essentially the same as that in \cite{subsampling99} immediately preceding Theorem 1 (of that paper), with the only difference
being that here we use a moment functional instead of a quantile functional.\footnote{For our problem, we found the moment functional to produce 
somewhat stabler estimates of $\zeta$ as compared to the quantiles.} 

\item[(4)] Estimate the $\alpha/2$'th and $(1-\alpha/2)$'th quantiles of $L$, say $q^{\star}_{\alpha/2}$ and $q^{\star}_{(1-\alpha)/2}$, from the empirical distribution of $n_j^{(1+2\hat{\zeta})/3}(\hat{\theta}_j^{\star} - \hat{\theta}_n)$ conditional on the data (either for $j=1$ or 2). This can be done by drawing a new set of subsamples of size $n_j$ from the original data.  

\item[(5)]  An approximate level $1-\alpha$ CI for $\theta^n$ is 
$[\hat{\theta}_n - q^{\star}_{(1-\alpha)/2}\,n^{-(1+2\hat{\zeta})/3}\;,\;\hat{\theta}_n - q^{\star}_{\alpha/2}\,n^{-(1+2\hat{\zeta})/3}]\,.$
\end{itemize} 
As this is a simple adaptation of an established procedure, we have not presented extensive simulation studies in the paper. However, we present results and figures from  limited simulation studies to provide a feel for the procedure. Data are generated from the model $Y=f_n(X) + \epsilon$, where 
\[ f_n(x) = \frac{\exp(n^{\gamma}(x - 0.5))}{1 + \exp(n^{\gamma}(x - 0.5))} \,,\] 
the random noise $\epsilon$ follows a normal distribution with mean zero and standard deviation $1.8$ and the covariate $X$ follows a uniform distribution on $[0,1]$. The 
sample size $n$ is taken to be 2000. Three values of $\gamma: 0.8, 1$ and $1.2$ are considered to account for each of the three regimes.  For demonstration purposes, we fit a one-parameter change-point model using $\beta_l^0=0$ and $\beta_u^0=1$ (as at the beginning of Section 4) and only estimate the change-point parameter. Given a dataset of size $2000$ from the above model with parameter $\gamma$, to estimate
$\hat{\zeta}$, we consider subsample sizes $n_1=100$ and $n_2=200$, draw $r=1000$ subsamples for each subsample size and then apply the formula in Step (3) above. The process is repeated for 200 datasets, resulting in 200 estimates $\{\hat{\zeta}_i\}_{i=1}^{200}$ and their median $\hat \zeta_m$ is chosen as the final estimate of $\zeta$. For the three settings, the estimated values are 0.76, $1$ and $1$ respectively, the corresponding true $\zeta$'s being $0.87, 1$ and $1$. 

Figure 1 presents QQ plots of the empirical distribution of $n^{\zeta}(\hat{\theta}_n - \theta^n)$\footnote{For the one parameter model $\theta^n$ is referred to as 
$\tilde{\theta}_n$ in Section 4.} (based on $1000$ independent datasets) versus the empirical subsampling distribution $\tilde{n}^{\hat{\zeta}_m}\,(\hat{\theta}^{\star} - \hat{\theta}_n)$ based on 1000 subsamples of size $\tilde{n} = 100$. For each of the three regimes, the plots show an approximate alignment with the $y=x$ line as would be expected. The empirical coverage probabilities for $\theta^n$ using subsampling based 95\% CIs  and the $\hat{\zeta}_m$ values from the previous paragraph (see the formula for the CI in Step (5) above) on 200 new data-sets (with $n = 2000$) are found to be $96.5\%$, $93.5\%$ and $93\%$ respectively for the three regimes. 

To demonstrate the performance of the \emph{fully adaptive} 5--step procedure described earlier in the section, we discuss results from a second simulation experiment from the  model above with $\gamma=0.5,~0.75,~1$ and $2$ and sample size $n=1500$. Given a dataset of size $1500$ from the regression model with parameter $\gamma$, in the first step, we resample the data $5000$ times with subsample sizes $n_1 = 100$ and $n_2 = 200$ respectively. For each subsample size, we evenly split these $5000$ subsamples into $10$ groups. We then compute 10 estimates of $\zeta$ via the formula in the Step 3 of the above 5--step procedure, each estimate using 500 (this is the $r$ from the general description of the procedure) subsamples of size 100 and 500 of size 200 and prescribe the median of these estimates as the value of $\hat{\zeta}$ to be used for the construction of the confidence interval for $\theta^n$. (Using the median provides additional stability to the estimation of $\zeta$.) The confidence interval construction (Steps 4 and 5) uses $500$ additional subsamples from the same dataset with subsample size $100$. The empirical confidence intervals for the four scenarios based on 200 independent datasets (average lengths of the 200 CIs in brackets) are $94\% (.129), ~96.5\% (.086),~97\% (.086) $ and $95\%(.081)$ respectively, providing numerical evidence of the proposed adaptive inference procedure. Note that the coverages reported in this paragraph are more realistic than the ones in the previous paragraph, since the CI for each dataset is based on an estimate of $\zeta$ \emph{computed from the same dataset} as is always the case in a real application. 

The adaptive procedure is computationally fairly intensive owing to the estimation of the convergence rate for each sample. Also note that the subsampling procedure,  by the very nature of it, involves tuning parameters (the subsample sizes) and this typically plays an important role in the reliability of the results (see the discussion towards the end of Section 4 in \cite{subsampling99}). Further investigations to fine tune and objectify the selection of subsample size in the context of subsampling with unknown convergence rates in general problems, and more specifically, in the change point problem we study in this paper would be very interesting but fall outside the scope of the current paper. 

%

\begin{figure}[h]\label{f1}
\begin{center}
\includegraphics[height=4in, width=6in]{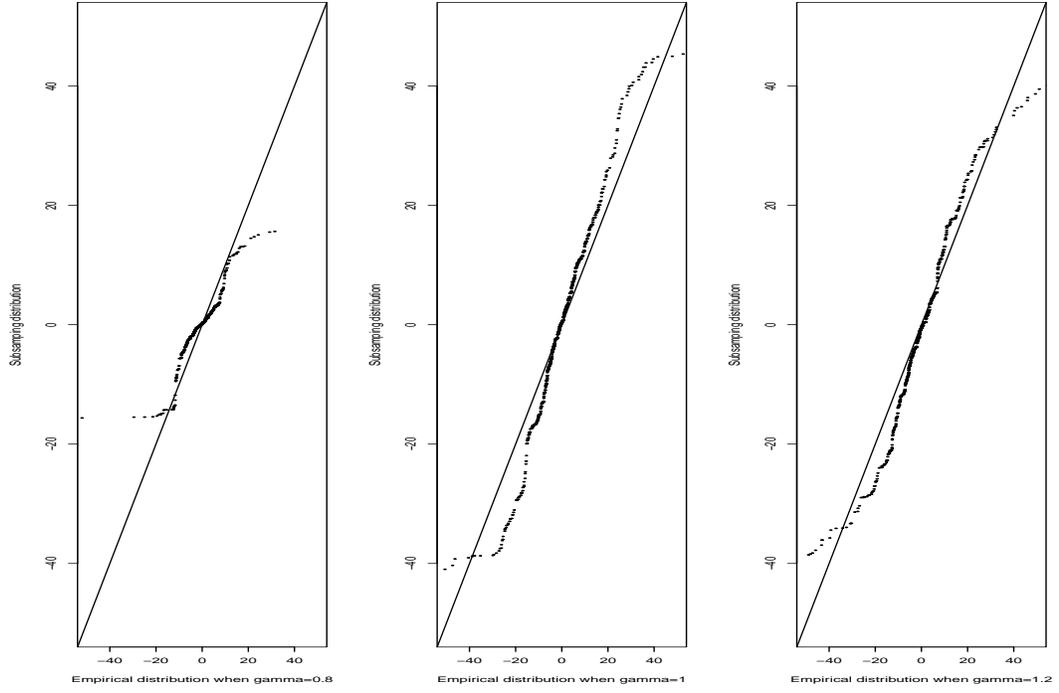}
\caption{The QQ-plots of subsampling empirical distributions versus the three empirical distributions for  sample size $n=2000$ and subsample size $100$. The true convergence rates are $n^{0.867}$, $n$ and $n$ respectively.  The straight line is a 45 degree line through the origin.}
\end{center}
\end{figure}




\section{Simulation studies}
In this section, we provide detailed empirical evidence of the theoretical results. Our framework stipulates a sequence of models changing with $n$ and converging to a limiting stump model with $\alpha_n$ regulating how fast the regression functions converge to a stump. We view a changing sequence of models with a given $\alpha_n$ as an asymptotic framework within which to couch \emph{a given fixed regression problem}: $Y = v_0(X) + \epsilon$ with $n$ data points available, and ask the question: which asymptotic regime: ``slow'' (i.e $\alpha_n < n$), ``intermediate'' ($\alpha_n = n$) or ``fast'' ($\alpha_n > n$) provides the best description of the behavior of $\hat{\theta}_n$, the least squares estimate of the population parameter $\theta_{0,v}$. Note that the population parameters $(\beta_{l,v}, \beta_{u,v}, \theta_{0,v})$ are given by the minimizer of $E(Y - \beta_{l}\,1(X \leq \theta) - 
\beta_u\,1(X > \theta))^2$ over all $(\beta_l, \beta_u, \theta)$, the expectation being taken with respect to the joint distribution of $(Y,X)$ in the above regression model. 

In the interests of a clean exposition, we restrict ourselves to three specific regimes: the one corresponding to $\alpha_n = c_0$ for a constant $c_0$ i.e. the fixed function set-up of Banerjee and McKeague (2007) which yields a Chernoff limit, the one with $\alpha_n = n$ that gives the intermediate distribution (Theorem 3.4)  and the last with $\alpha_n> n$, which produces the compound Poisson process limit which also arises when the true regression model is a fixed change-point model (Theorem 3.5). The case $\alpha_n = c_0$ should be viewed as a representative of the slow regimes corresponding to $\alpha_n < n$: recall that all slow regimes lead to a multiple of Chernoff's distribution, albeit with different convergence rates.

We generate data from the model $Y=v_0(X) + \epsilon$, where 
\[ v_0(x) = \frac{\exp(M(x - 0.5))}{1 + \exp(M(x - 0.5))} \,,\] 
the random noise $\epsilon$ follows a normal distribution with mean zero and standard deviation $0.6$, and $M$ is a constant that we vary for different simulation settings as explained below.  The covariate $X$ follows a piecewise uniform distribution on $[0,1]$ and is symmetric about $0.5$. To demonstrate the effect of the density $p_{X}$ on the limiting distribution, we consider two scenarios for generating the covariate $X$. In scenario 1, the density is $8$ on $[0.45, 0.55]$, whilst in scenario 2, the density is $4$ on  $[0.4, 0.6]$. For scenario 1, we consider four different values of $M: 35, 60, 100, 1000$, and for scenario 2, $M$ assumes values $20, 60, 100, 1000$. Larger values of $M$ produce steeper sigmoidal curves which are closer to a change point model than smaller values. For each combination of $p_X$ and $M$ (leading to 8 settings), we generate data for sample sizes $n = 50, 100, 1000, 4000$; for each sample, we generate $500$ data sets (replicates) to get the empirical distribution of $\hat{\theta}^n$. For the limiting distribution based on the fixed function setting, the normalized least squares estimate $n^{1/3}(\hat{\theta}^n - \theta^n)$ is calibrated against the appropriate Chernoff limit; see Theorem 2.1 of \citep{BM:07}. For the limiting distribution based on $\alpha_n=n$, we calibrate $n\,(\hat{\theta}^n - \theta^n)$ against the quantiles of the argmin of $\Lambda$ in Theorem 3.4; whilst, for the third case, we calibrate $n\,(\hat{\theta}^n - \theta^n)$ against the quantiles of the argmin of $\Lambda_1$ in Theorem 3.5. Note that in our asymptotic framework, $\theta^n$ is the population minimizer of the change--point parameter at stage $n$, and since we have a fixed regression model $v_0$ in our set-up, this is identically equal to the parameter $\theta_{0,v}$. 

Simulating from the limit distribution in Theorem 2.1 of  \citep{BM:07} -- the slow regime -- requires the $(\beta_l^0, \beta_u^0, d^0)$ appearing in that result. For us, these are simply the population parameters $(\beta_{l,v}, \beta_{u,v}, \theta_{0,v})$ (which depend on $M$ and $p_X$). To simulate the theoretical limiting distribution based on Theorem 3.4 (intermediate regime), we write our fixed function $v_0$ as $v_0(x) = f(n\,(x- 0.5))$ (so as to obtain the representation $Y = f(n(X-0.5)) + \epsilon$ based on which the limit is derived); here, $f$, of course, becomes dependent on $n$: $f_n(t) = \frac{\exp(Mt/n)}{1+\exp(Mt/n)}$. The quantity $\xi_0 = f_n^{-1}((\beta_l^0 + \beta_u^0)/2$, needed to generate $\Lambda$ in Theorem 3.4  is $0.5$, since $\beta_u^0$ and $\beta_l^0$, the limits of $f_n$ as $t$ goes to $\infty$ and $-\infty$ respectively are $1$ and $0$. To simulate the limit distribution based on Theorem 3.5, we require the value $(\beta_u^0 - \beta_l^0)/2$ in that theorem, and here $\beta_u^0 = 1$ and $\beta_l^0 = 0$ are the levels of the limiting change-point model. 

The QQ-plots of the empirical distribution of the normalized least squares estimate based on 500 replicates against that of a sample drawn from the limiting distribution for each regime (the size of the sample from the limiting approximation is 2000 in every simulation setting) are presented in a series of figures: two of these corresponding to scenario 1, $M = 1000$ and scenario 2, $M = 1000$ are presented in the main paper and the rest are included in \cite{Supple}. The general pattern is fairly clear: the fast regime provides better approximations at smaller sample sizes than at larger ones, the slow regime improves at higher sample sizes, and the intermediate regime is much more robust to the sample size, though at high sample sizes ($n = 1000, 4000$) the approximation provided by it starts breaking down (see, for example, the behavior of the intermediate regime for smaller values of $M$). The fast regime generally completely breaks down at high sample sizes and for smaller values of $M$ (20, 35, 60), which correspond to curves that are farther from a change-point model, tends to behave poorly even for small samples. While the slow regime improves for larger samples, it sometimes provides a decent approximation at smaller samples as well (again, see some of the plots in \cite{Supple}). 

The general pattern can be explained by noting that at small sample sizes the data is typically not adequate to discover the features of the underlying sigmoidal curve; especially for a steep curve (for example $M = 1000$ as presented in the paper), at small $n$, the data only `sees' the change-point type feature, and therefore an approximation using the fast regime (also the regime for a true change-point model) performs better. For large $n$, the data is able to `pick out' the overall pattern of the continuous curve quite well and consequently, the setting of \citep{BM:07} which deals with fitting a change-point working model to a smooth fixed regression function is apt. The intermediate setting or the boundary case strikes a balance between these two approximations as it uses some features of the underlying regression curve but on the other hand not as local features as the ones used by the asymptotics in \citep{BM:07}. Hence, it provides an approximation that adapts much better to changes in sample size. This is consistent with the fact that the family of boundary distributions can transition to either of the two extreme limits, as shown in Section 4.

\begin{figure}[h]\label{f1}
\begin{center}
\includegraphics[height=7in, width=6in]{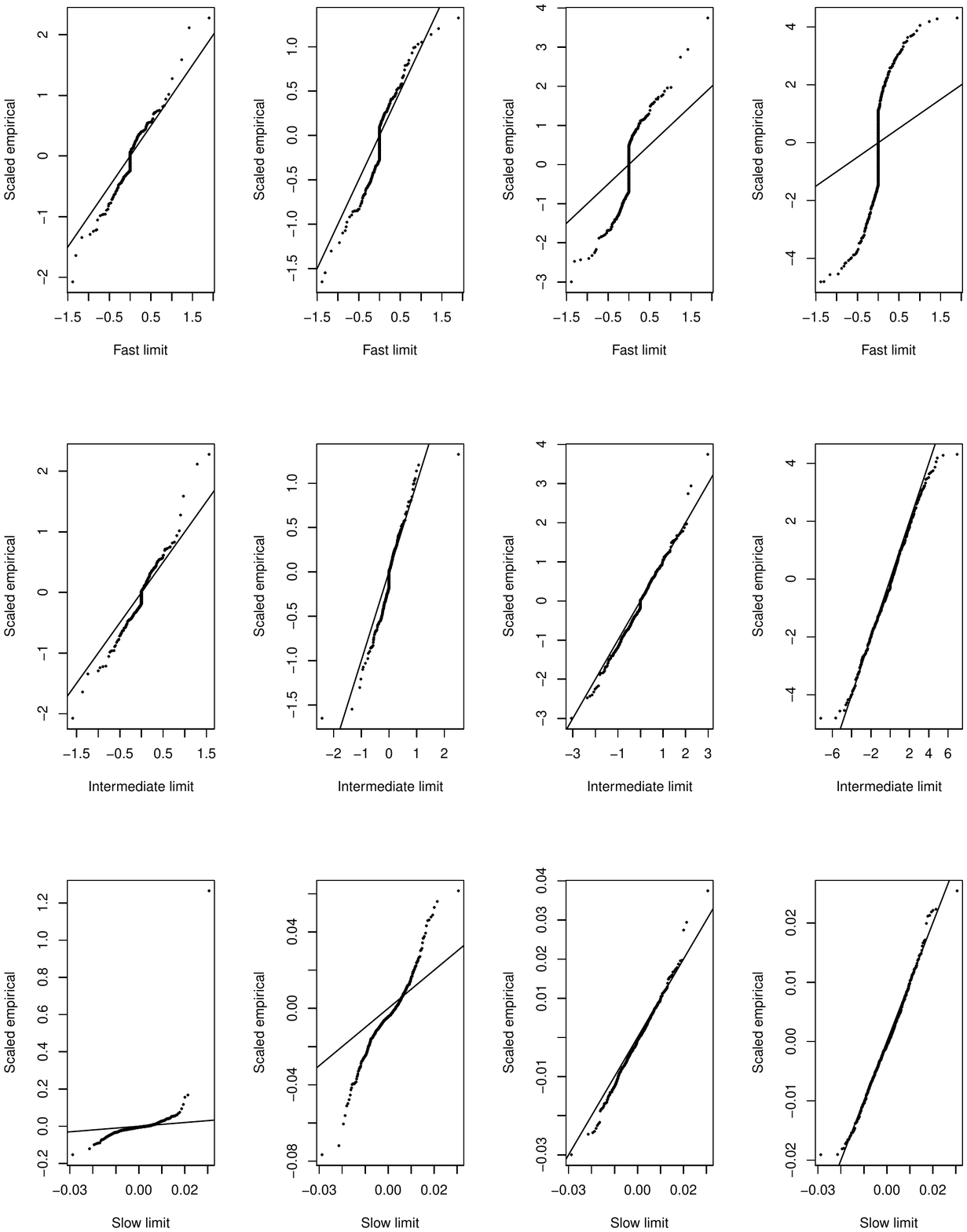}
\vspace{-1cm}
\caption{The QQ-plots of empirical distributions versus the three theoretical limits for $p_X(\theta_0)=8$ and $M=1000$. First row and second row:  empirical distribution of $n(\hat \theta^n - \theta^n)$ vs the fast limit and the intermediate limit respectively. Third row: empirical distribution of $n^{1/3}(\hat \theta^n - \theta^n)$ vs the slow limit. The four columns represent $n=50, ~100,~1000,~$ and $4000$ respectively.  The straight line is a 45 degree line through the origin.}

\end{center}
\end{figure}

\begin{figure}[h]\label{f1}
\begin{center}
\includegraphics[height=7in, width=6in]{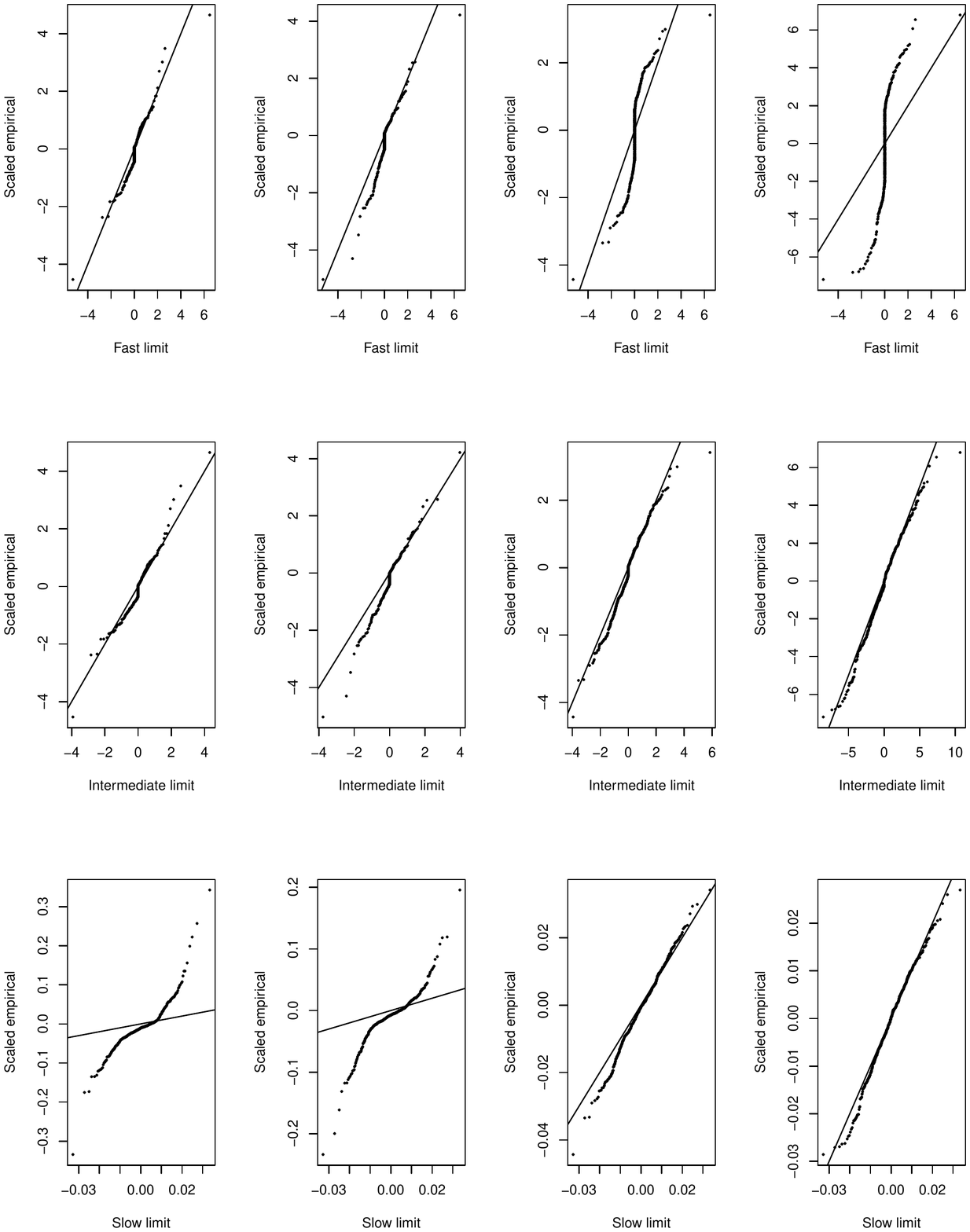}
\caption{The QQ-plots of empirical distributions versus the three theoretical limits for $p_X(\theta_0)=4$ and $M=1000$. The rest of the captions and notations are the same as those in Figure 1.}

\end{center}
\end{figure}

\section{Discussion}
In this paper, we have studied the asymptotic behavior of change-point models under a wide range of model mis-specification. We end with a discussion of some important aspects of our work and some related problems. 

{\bf Analogy to classical parametric models:} Viewing the $f_n$'s as a sequence of local
alternatives to the limiting null model: the stump function $f_0(x; \psi_0)=\beta_l^0 1(x \le \theta^0) + \beta_u^0 1(x> \theta^0)$, the phenomena studied in this paper are qualitatively identical to what transpires with the MLE in regular parametric models under a sequence of local alternatives. 

So, consider such a model $\{p(x; \eta)\}$ with the p-dimensional parameter $\eta$ and let $X_1, X_2, \ldots X_n$  denote i.i.d. observations. Let $\hat \eta$ denote the MLE for $\eta$. We aim to test the null hypothesis $ \eta = \eta_0$. It is well known that under the null, $\sqrt{n}(\hat \eta - \eta_0)$ follows an asymptotic normal distribution $N(0, I(\eta_0)^{-1})$, where $I(\eta_0)$ is the information matrix for $\eta$. With alternatives converging faster than $\sqrt{n}$, say $\eta_n = \eta_0 + h\,n^{-\gamma}$ for $\gamma > 1/2$, the limit of $\sqrt{n}(\hat{\eta}_n - \eta_0)$ continues to be identical to that under the null. With alternatives converging at a slower than the `regular' $\sqrt{n}$ rate, i.e. when $\gamma < 1/2$, the limit distribution of $\sqrt{n}\,(\hat{\eta}_n - \eta_0)$ is no longer tight, since the bias term $\sqrt{n}(\eta_n - \eta_0)$ drifts to $\infty$. In the change-point problem, $f_0(x,\psi_0)$, of course, plays the role of $\eta_0$, the convergence rate $n$, which is the natural convergence rate of the least-squares estimate of $\theta^0$ under the null model, plays the role of $\sqrt{n}$, the $f_n$'s take on the role of $\eta_n$, while $\alpha_n$ becomes the analogue of $n^{-\gamma}$. As noted in the discussion before the statement of Theorem \ref{the-weak-fastn}, for $\alpha_n$ going to $\infty$ faster than $n$ (corresponding in the classical case to $\gamma > 1/2$), the asymptotic distribution of $\hat{\theta}^n$ in our problem is identical to that under the null model $f_0$. When $\alpha_n = o(n)$ (corresponding in the classical case to $\gamma < 1/2$), Theorem \ref{the-weak-smalln} in conjunction with Theorem \ref{the-det-rate} tells us that $n^{1/3}\,\alpha_n^{2/3}(\hat{\theta}^n - \theta_0)$ does not have a tight limit, since the bias term $n^{1/3}\,\alpha_n^{2/3}(\theta^n - \theta^0)$ goes to $\infty$. 

It remains to compare the cases where the alternative approaches the null at the natural convergence rate. In the classical scenario, this corresponds to $\gamma = 1/2$ and produces a tight distribution in the limit, namely, 
$N(h, I(\eta_0)^{-1})$ for $\sqrt{n}(\hat{\eta}_n - \eta_0)$; thus, the direction of approach of the local alternatives figures in the limit. In the change--point scenario, the analogous situation is $\alpha_n = n$, and as Theorem \ref{the-weak-n} shows, now the distribution of $n(\hat{\theta}^n - \theta^0)$ converges to a tight limit which depends upon $f$, which can be interpreted as the `direction' in which the smooth $f_n$'s approach the stump $f_0$. One important difference between the classical and the change--point scenario is, of course, the differing convergence rates: the $\alpha_n$ parameter influences the rate at which $\hat{\theta}^n$ approaches 
$\theta^n$ in the change--point model, but the $\gamma$ parameter in the classical scenario does not influence the convergence rate: in fact, $\sqrt{n}(\hat{\eta}_n - \eta_n)$ is $O_p(1)$ in all situations. 
\newline
\newline
{\bf An alternative approach for inference:} 
An alternative approach for inference in this problem, kindly brought to our attention by a referee, relies on 
smoothed least squares estimation along the lines pursued in the papers \cite{firstseopaper} and \cite{secondseopaper}. \cite{firstseopaper} studies a linear regression model with regime switching, where the form of the linear regression depends on whether a particular subset of covariates lies above or below a hyperplane (whose parameters are also unknown). This can be thought of as a `change-plane' problem. To avoid the non-standard distributions that would come into play under a regular least squares approach, the authors replace the indicator function appearing in the least squares criterion by a smooth integrated kernel function (analogous to a distribution function) in the spirit of the smoothed maximum score estimator of \cite{horowitz1992smoothed}. The corresponding smoothed least squares estimators -- even those of the hyperplane parameters -- are seen to be asymptotically normal under appropriate conditions on the model and the bandwidth used for the integrated kernel function. Asymptotic normality makes inference more tractable though the rate of convergence is somewhat compromised and can be at most $n^{3/4}$ (up to a logarithmic factor), slower than the $n$--rate of convergence attained by the regular least squares estimators. While the set-up of \cite{firstseopaper} works under the assumption that the threshold model defined by the hyperplane is true, \cite{secondseopaper} (Section 4) explores the behavior of the smoothed least squares estimate under mis-specification in the spirit of our paper and establishes asymptotic normality (Theorem 4), with the rate again depending upon the bandwidth used. 
These investigations suggest the possibility that using a smoothed least squares approach in our diminishing mis-specification problem could also lead to asymptotic normality, avoiding the non-standard distributions that now come into play, at the expense of somewhat reduced convergence rates. 
\newline
\newline
{\bf Other potential extensions and connections:} A natural question is the extension of this approach to multiple change points , i.e. a situation where the limit of the converging (smooth) models is a piecewise constant function
with multiple jumps. It is clear that the properties of the underlying $f$ (i.e. Assumptions A through C) which were used to manufacture the converging models would now need to change. Recall that in this paper, the regression function at stage $n$, $f_n(x) = f(\alpha_n(x - \theta_0))$ and as $\alpha_n$ goes to $\infty$, $f_n$ must necessarily converge to a piecewise constant function with a single jump. For example, to take into account the situation where the limiting function is of the form 
$\alpha_0\,1(x < \theta_{10}) + \beta_0\,1(\theta_{10} \leq x < \theta_{20}) + \gamma_0\,1(\theta_{20} \leq x)$, one possibility for a converging smooth function $\tilde{f}_n$ could be: 
\[f_n(x) = f_1(\alpha_n(x - \theta_{10})) + f_2(\alpha_n(x - \theta_{20})) \,,\]
for monotone functions $f_1, f_2$ and appropriate conditions on their limit values at $-\infty$ and $\infty$. Note, moreover, that the $\alpha_n$ sitting within the $f_2$ could be replaced by a different rate parameter ($\beta_n \ne \alpha_n$)
 going to $\infinity$. Thus, the multiple change point problem throws up a number of different challenges which are outside the scope of this paper. 

In conclusion we would like to note an interesting connection of our results to \cite{F2007}, also pointed out by a referee. Section 3.2 of \cite{F2007} considers approximating functions in different smoothness classes using the Unbalanced Haar transform as basis vectors for the approximating class. These basis functions are piecewise constant by construction and are therefore expected to provide more precise approximations to underlying functions that are structurally similar. Indeed, the result of Theorem 3.1 in that paper shows that when $f$ is in the class of functions of bounded variation the expected IMSE of the Haar transform based estimate attains a rate of $n^{-2/3}$ up to a logarithmic factor that involves the sample size as well as certain features of the approximation basis. On the other hand, when $f$ is in $S[0,1]$, the space of piecewise constant functions with finitely many jumps, the rate improves to $n^{-1}$, again up to logarithmic terms and the number of jumps of the function. 
In our work, the approximating function is a piecewise constant function with a single jump (a stump) and the underlying function a smooth function that can be considered close to a (limiting) stump with a jump at $\theta_0$. The degree of closeness is measured by the parameter $\alpha_n$. Our results show that for larger values of 
$\alpha_n$, which correspond to the underlying function behaving more like a stump, the rate of convergence of $\hat{\theta}_n$ is faster: for $\alpha_n$ at least as large as $n$, $|\hat{\theta}_n - \theta_0| = O_p(1/n)$ 
whereas for $\alpha_n = o(n)$, $|\hat{\theta}_n - \theta_0| = O_p(\alpha_n^{-1})$ and therefore slower than the other case. 
We note, of course, that in contrast to \cite{F2007} where a global measure of error is considered, our results are formulated in terms of the convergence of the estimated jump parameter alone.






\bigskip


\bibliographystyle{ims}
\bibliography{mischange}
\end{document}